%% file: CrossedModules.tex
\documentclass[11pt,a4paper,reqno]{amsart}

\usepackage{amssymb}
\usepackage{amsmath}

\usepackage[matrix,arrow,tips,curve]{xy}
\SelectTips{cm}{11}

\usepackage{rotating}                      
\newcommand{\rotateninety}[1]{\begin{rotate}{90}$#1$\end{rotate}}

\newcommand{\formularemark}[1]{\text{\scriptsize #1\quad}}

\newcommand{\mysection}[1]{\vspace{1em}\section[#1]{\uppercase{#1}}}
\newcommand{\mysectionstar}[1]{\vspace{1em}\section*[#1]{\uppercase{#1}}}

\theoremstyle{plain}
\newtheorem*{proposition}{Proposition}
\newtheorem*{propositiondefinition}{Proposition and Definition}
\newtheorem*{theorem}{Theorem}
\newtheorem*{corollary}{Corollary}
\theoremstyle{definition}
\newtheorem*{definition}{Definition}
\newtheorem*{fact}{Definition and Fact}
\theoremstyle{remark}
\newtheorem{example}{Example}
\newtheorem*{example*}{Example}
\newtheorem*{remark}{Remark}

\newcommand{\ti}[1]{\mbox{\b{$#1$}}}
\newcommand{\Aut}[1]{\ensuremath{\operatorname{Aut}#1}}
\newcommand{\Ad}[1]{\ensuremath{\operatorname{Ad}_{#1}}}

\newcommand{\twogroup}[1]{\ensuremath{\mathcal{#1}}}

\newcommand{\GG}{\twogroup{G}}
\newcommand{\CC}{\twogroup{C}}
\newcommand{\twoGrp}{\textit{2-Grp}}
\newcommand{\XMod}{\textit{X-Mod}}
\newcommand{\Cat}{\textit{Cat}}
\newcommand{\Fun}{\textit{Fun}}
\newcommand{\coker}{\ensuremath{\operatorname{coker}}}
\newcommand{\id}{\ensuremath{\operatorname{id}}}
\renewcommand{\phi}{\varphi}
\renewcommand{\theta}{\vartheta}
\newcommand{\restrictedto}[1]{_{|#1}}

\input{sspbibstyle}

\usepackage[bookmarksopen]{hyperref}
\pdfcatalog{
	/Pagemode	/UseOutlines
	/OpenAction 	/fitbh
}

\begin{document}

\renewcommand{\abstractname}{ABSTRACT}
\begin{abstract} 
The 2-categories of strict 2-groups and crossed modules are introduced  and their 2-equivalence is made explicit.
\end{abstract}

\title{Strict 2-Groups are Crossed Modules}
\author{Sven-S. Porst}
\thanks{\href{mailto:ssp@uni-math.gwdg.de?subject=Crossed\%20Modules}{ssp@uni-math.gwdg.de}\\
\hspace*{1em} Mathematisches Institut $\cdot$ Bunsenstra\ss e  3--5 $\cdot$ D--37073 G\"{o}ttingen}
\date{}
\maketitle

\thispagestyle{empty}

%
%
\mysectionstar{Introduction}
\label{introduction}
\noindent
Crossed modules and strict 2-groups 
have been known since the middle of the twentieth century. Their recent renaissance in the area of `higher' bundles and gauge theory brings them to then attention of differential geometers and mathematical physicists who may not be as familiar with the underlying concepts as homotopy or category theorists are.

This text aims to define crossed modules and strict 2-groups in a concrete way. Furthermore the frequently quoted 2-equivalence of the 2-categories of crossed modules and strict 2-groups is shown in full detail. In particular, the text includes a solution to the very last exercise in \cite{borceux-hca1} and the argument found in Forrester-Barker's helpful paper \cite{math.CT/0212065}. 
With the \emph{2}-equivalence being established, the result goes a bit beyond those and provides the omitted proof to Theorem 2 of \cite{brownspencer-ggcmatfg}. Category theorists will find the hands-on computations too laborious and may prefer the more abstract proof of the same fact in Theorem 5.13 of \cite{fiore-papdc}.

The computations needed to do this are technically simple and occasionally tedious. The main effort is keeping in mind the kinds of objects and morphisms being used. Many brackets and multiplication signs are omitted from notation for the sake of legibility. However, the formul\ae's meaning should still remain clear in the context of the elements and letters used.

An introduction to internal categories and 2-categories is given in section~\ref{categories}. It highlights the parallels between the definition of a category and that of a 2-category -- going along the lines of chapter~7 in \cite{borceux-hca1}. Another convenient introduction  is contained in chapter~XII of  \cite{maclane-cftwm} while an early review of the subject is given in \cite{kellystreet-roteo2c}.

Strict 2-groups are introduced in section~\ref{2groups} as internal categories in the category of groups. It follows immediately from this that the composition in a 2-group is uniquely determined by the rest of the structure. Crossed modules along with their morphisms and 2-morphisms are introduced in section~\ref{crossedmodules}. They come from homotopy theory and have been around for well over half a century \cite{whitehead-noapp,whitehead-ch2}; An account of that can be found in \cite{norrie-aaaocm}. 2-morphisms of crossed modules were already defined as homotopies in \cite{cockcroft-othos}.

In sections~\ref{2groupstomodules} and~\ref{modulesto2groups} 2-functors from strict 2-groups to crossed modules and vice versa are given. The construction used on the object level can also be found in section XII.8 of \cite{maclane-cftwm} and in more detail in \cite{math.CT/0212065}. Section~\ref{roundtrip} finally shows that these 2-functors establish a 2-equivalence.

%
%
\mysection{Categories, internal categories, 2-categories}
\label{categories}
\noindent
We begin by recalling basic definitions from category theory and introducing internal categories as well as 2-categories. 

\begin{definition}
A \emph{category} $\GG$ consists of a \emph{class of objects} $G_0$ and a \emph{set of morphisms} $G_1(A,B)$ for all pairs $(A,B)$ of objects.  For any triple of  objects $(A, B, C)$ there is a composition map
\[
\circ: G_1(B, C) \times G_1(A,B) \longrightarrow G_1(A,C)
\]
and for each object $A$ there is an identity morphism
\[
i A \in G_1(A,A)\,.
\]
The composition is required to be associative and the identity morphisms have to satisfy $i B \circ f = f = f \circ i A$ for all morphisms $f\in G_1(A, B)$.
\end{definition}

We denote by $G_1$ the class of all morphisms in $\GG$ and write the
\emph{source}, \emph{target} and  \emph{identity} maps
\[
s: G_1 \longrightarrow G_0
\qquad\qquad
t: G_1 \longrightarrow G_0
\qquad\qquad
i: G_0 \longrightarrow G_1
\]
We use the pullback $G_1 \!\ _{s}\!\times_{t} G_1$ to define the \emph{composition map} globally
\[
\circ: G_1 \!\ _{s}\!\times_{t} G_1 \rightarrow G_1\,.
\]
In this notation the conditions for $\GG$ being a category can be written in terms of diagrams only. They are equivalent to commutativity of the following diagrams:
\begin{gather}
\label{categoryconditions}
\xymatrix{
& G_0 
	\ar@/_.0001mm/[dl]_{\id_{G_0}}
	\ar[d]^i
	\ar@/_.0001mm/[dr]^{\id_{G_0}} 
& \\
G_0
& G_1
	\ar[l]^s
	\ar[r]_t
& G_0	
}
\qquad
\xymatrix@C=10mm{
G_1\times G_0 
	\ar[d]_{\id_{G_1}\times i}
& G_1
	\ar[l]_<>(.5){\id_{G_1}\times s}
	\ar[d]^{\id_{G_1}}
	\ar[r]^<>(.5){t \times \id_{G_1}}
& G_0 \times G_1 
	\ar[d]^{i\times \id_{G_1}}
\\
G_1 \!\ _{s}\!\times_{t} G_1
	\ar[r]_<>(.5)\circ
& G_1
& G_1 \!\ _{s}\!\times_{t} G_1
	\ar[l]^<>(.4)\circ
}
\\
\xymatrix{
G_1 
	\ar[d]_s
& G_1 \!\ _{s}\!\times_{t} G_1
	\ar[l]_<>(.4){\pi_2}
	\ar[d]^\circ
	\ar[r]^<>(.4){\pi_1}
& G_1
	\ar[d]^t
\\
G_0
& G_1
	\ar[l]^s
	\ar[r]_t
& G_0
}
\qquad
\xymatrix@C=12mm{
G_1 \!\ _{s}\!\times_{t} G_1 \!\ _{s}\!\times_{t} G_1
	\ar[d]_{\id_{G_1} \times \circ}
	\ar[r]^<>(.5){\circ \times \id_{G_1}}
& G_1 \!\ _{s}\!\times_{t} G_1 
	\ar[d]^\circ
\\
G_1 \!\ _{s}\!\times_{t} G_1 
	\ar[r]_<>(.5)\circ
& G_1
}
\notag
\end{gather}

\begin{definition}
Given two categories $\GG$ and $\GG'$, a \emph{functor} $F: \GG \rightarrow \GG'$ consists of two maps
\[
F_0 : G_0 \longrightarrow G_0' 
\qquad \quad \text{and} \qquad \qquad
F_1: G_1 \longrightarrow G_1'
\]
which are compatible with the source, target, identity as well as composition maps, meaning that the following diagrams commute:
\begin{gather}
\label{functorconditions}
\xymatrix@C=7mm@R=7mm{
G_0 	
	\ar[d]_{F_0}
& G_1	
	\ar[l]_s
	\ar[d]^{F_1}
	\ar[r]^t
& G_0
	\ar[d]^{F_0}
\\
G_0'
& G_1' 
	\ar[l]^{s'}
	\ar[r]_{t'}
& G_0'
}
\quad\;\;\;
\xymatrix@C=7mm@R=7mm{
G_0 
	\ar[d]_{F_0}
	\ar[r]^i
& G_1
	\ar[d]^{F_1}
\\
G_0'
	\ar[r]_{i'}
& G_1'
}
\quad\;\;\;
\xymatrix@C=7mm@R=7mm{
G_1 \!\ _{s}\!\times_{t} G_1
	\ar[d]_{F_1\times F_1}
	\ar[r]^<>(.5)\circ
& G_1
	\ar[d]^{F_1}
\\
G_1' \!\ _{s'}\!\times_{t'} G_1'
	\ar[r]_<>(.5){\circ'}
& G_1'
}
\end{gather}
\end{definition}
We shall omit the index of $F$ from notation whenever the context makes clear whether we are working on the level of objects or morphisms.

\begin{definition}
Given two functors $F, E: \GG \rightarrow \GG'$, a \emph{natural transformation} $\theta: F\Rightarrow E$ assigns to every object $A$ in $\GG$ a morphism $\theta A: F A \rightarrow E A$ in $\GG'$ \\[1mm]
\noindent
\begin{minipage}{55mm}
such that for each morphism\\
$f\in G_1(A,B)$ 
diagram~\eqref{naturaltransformation} commutes. 
If all $\theta A$ are isomorphisms, $\theta$ is
called a \emph{natural isomorphism}.

Diagram  \eqref{naturaltransformation} is a diagram in $\GG$ and its commutativity is equivalent to commutativity of diagram \eqref{naturaltransformation2} in the category of sets. This will be useful when generalising this definition to 2-categories in a moment.
\end{minipage}
\quad
\begin{minipage}[c]{66.5mm}
\begin{equation}
\stepcounter{equation}
\tag{\theequation a}
\label{naturaltransformation}
\xymatrix@C=20mm{
F A 	
	\ar[r]^<>(.5){F f} 
	\ar[d]_{\theta A}
& F B  
	\ar[d]^{\theta B} 
\\
E A 
	\ar[r]_<>(.5){E f }
& E B 
}
\end{equation}
\begin{equation}
\tag{\theequation b}
\label{naturaltransformation2}
\xymatrix@C=7mm{
\GG(A,B)
	\ar[r]^<>(0.5){F}
	\ar[d]_{E}
& \GG'(F\!A, FB)
	\ar[d]^{\theta B \circ} \\
\GG'(E A, E B)
	\ar[r]_<>(0.5){\circ \theta A}
& \GG'(F\! A, E B)
}
\end{equation}
\end{minipage}
\end{definition}

Again, the condition of a map $\theta: G_0 \rightarrow G_1'$ being a natural transformation can be expressed in the language of diagrams. It is equivalent to the commutativity of the following two: 
\begin{equation}
\label{ntconditions}
\xymatrix@C=15mm{
& G_0
	\ar@/.000001mm/[ld]_{F_0}
	\ar[d]^{\theta}
	\ar@/.000001mm/[rd]^{E_0}
& \\
G_0'
& G_1'
	\ar[l]_{s'}
	\ar[r]^{t'}
& G_0'
}
\qquad
\xymatrix{
G_1
	\ar[r]^<>(.5){(\theta t, F_1)}
	\ar[d]_{(E_1, \theta s)}
& G_1' \!\ _{s'}\!\times_{t'} G_1'
	\ar[d]^{\circ'}
\\
G_1' \!\ _{s'}\!\times_{t'} G_1'
	\ar[r]_<>(.5){\circ'}
& G_1'
}
\end{equation}

\begin{definition}
Two categories $\GG, \GG'$ are called \emph{equivalent} if there are functors
\[
F: \GG \longrightarrow \GG'
\qquad \qquad \text{and} \qquad\qquad
E: \GG' \longrightarrow \GG
\]
such that there are natural isomorphisms 
\[
\theta: E \circ F \stackrel{\!\sim}{\Longrightarrow} \id_{\GG}
\qquad \qquad \text{and} \qquad \qquad
\sigma: F \circ E \stackrel{\!\sim}{\Longrightarrow} \id_{\GG'}\,.
\]
\end{definition}

\begin{definition}
Given a  category $\CC$ with finite limits, we say that the tuple of objects and morphisms in $\CC$ $\GG = (G_0, G_1, s, t, i, \circ)$ is an \emph{internal category in $\CC$}
if the four diagrams in \eqref{categoryconditions} commute.

An \emph{internal functor in $\CC$} between internal categories $\GG$ and $\GG'$ is a  pair of $\CC$-morphisms $(F_0: G_0 \rightarrow G_0', \; F_1: G_1 \rightarrow G_1')$ such that the diagrams in \eqref{functorconditions} commute. 
An \emph{internal natural transformation in $\CC$} between internal functors $F, F': \GG \rightarrow \GG'$ is a morphism $\theta: G_0 \rightarrow G_1'$ in $\CC$ making the diagrams in
\eqref{ntconditions} commute.
\end{definition}

In general we have to express the conditions for an internal category in the language of morphisms and commutative diagrams. But given a faithful forgetful functor $F: \CC \rightarrow \mathrm{Sets}$ which preserves and reflects pullbacks, the tuple $(G_0, G_1, s,t, i, \circ)$ is an internal category in $\CC$ if and only if $(FG_0, FG_1, Fs, Ft, Fi, F\circ)$ is a small category. That is,  if there is a suitable forgetful functor to the category of sets, we can work with an internal category as if it were a category. The forgetful functor from the category of groups to the category of sets is an example for this which will implicitly be used in the following sections.

\begin{definition}
A \emph{2-category} $\GG$ consists of  a class of 
objects $G_0$ and for any pair of  objects $(A, B)$ a small \emph{category of morphisms} $\GG(A,B)$ -- with objects $G_1(A,B)$ and morphisms $G_2(A,B)$ --, along with \emph{composition functors} 
\[
\circ: \GG(B,C) \times \GG(A,B) \longrightarrow \GG(A,C)
\]
for every triple  $(A, B, C)$ of  objects and \emph{identity functors} from the terminal category 
\[
iA: 1 \longrightarrow \GG(A,A)
\]
for all objects $A$ such that $\circ$ is associative 
and 
\[
i B \circ F = F = F \circ i A
\qquad\text{ as well as } \qquad
i B \circ \theta = \theta = \theta \circ i A
\] 
hold for all $F\in G_1(A,B)$ and $\theta \in G_2(A,B)$.

For all pairs of objects $(A,B)$ elements  of $G_1(A,B)$ are called \emph{morphisms} or \emph{1-cells} of $\GG$ and elements of $G_2(A,B)$ are called \emph{2-morphisms} or \emph{2-cells} of $\GG$. We write $G_1$ and $G_2$ for the classes of  all morphisms and 2-morphisms respectively.
 
There are two ways of composing 2-morphisms: using the composition $\bullet$ inside the categories $\GG(A,B)$, called \emph{vertical} composition, and using the morphism level of the  functor $\circ$, called \emph{horizontal} composition. 
\end{definition}

\begin{remark}
A 2-category $\GG$ gives rise to various category structures: the morphism category $\GG(A,B)$ for each pair $(A,B)$ of objects; the category of objects and morphisms between them $(G_0, G_1)$; and the category of objects and 2-morphisms between them $(G_0, G_2)$.
This means we have three different compositions already. In the following example our objects will be categories and thus each of them brings its own composition as well. 
\end{remark}

\begin{example*}
The standard example is the 2-category $\Cat$ of small categories, functors and natural transformations. For any pair of small categories $(\GG, \GG')$ we have the functor category $\Fun(\GG, \GG')$. Compositions inside the $\Fun(\GG, \GG')$  give the vertical composition of the 2-category. Its horizontal composition is given by
\begin{equation}
\label{horizontalcompositionnaturaltransformations}
\begin{array}{rcll}
\circ: & \Fun(\GG', \GG'') \times \Fun(\GG, \GG') &\longrightarrow &\Fun(\GG, \GG') \\
& (F', F) & \longmapsto  & F'\circ F\\
& (\theta': F' \Rightarrow E', \; \theta: F \Rightarrow E) & \longmapsto & 
[A \mapsto \theta'(E A) \, F' \!( \theta A )] \,,
\end{array}
\end{equation}
\begin{minipage}{0.52\textwidth}
that is by composition of functors on the object level; To understand the construction on the morphism level, consider the natural transformation diagram \eqref{naturaltransformation} for the natural transformation $\theta'$ applied to the morphism $\theta A: FA \rightarrow EA$ in $\GG'$.\\[-3.2mm]
\end{minipage}
\begin{minipage}{0.48\textwidth}
\[
\xymatrix@C=20mm{
F'FA \ar[r]^{F'(\theta A)}
	\ar@{-->}[dr]^{(\theta'\circ\theta)(A)}
	\ar[d]_{\theta'(FA)}
& F'EA \ar[d]^{\theta'(EA)} 
\\
E'FA	\ar[r]_{E'(\theta A)}
&	E'EA
}
\]
\end{minipage}

Its diagonal map is just the morphism defined above and commutativity of the diagram means that we can alternatively use $E'(\theta A) \, \theta'( FA)$. Once this is clear and identity functors have been defined as
\begin{equation}
\label{identitynaturaltransformation}
i\GG(*) = \id_{\GG}
\qquad \qquad
i\GG(\id_{*}) = [\,\id_{\id_\GG} : \id_\GG \Rightarrow \id_\GG \quad A \mapsto iA\,]
\end{equation}
it is quickly seen that all the required identities are satisfied.
\end{example*}

\begin{example*}
The previous example also holds for internal categories:  For a category $\CC$ with finite limits,  internal categories in $\CC$, along with functors in $\CC$ and natural transformations in $\CC$, form a 2-category. 
\end{example*}

\begin{example*}
Another classical example for a 2-category is given by topological spaces, continuous maps and homotopy classes of homotopies.
\end{example*}

\begin{definition}
Given two 2-categories $\GG, \GG'$, a \emph{2-functor} $F:\GG \rightarrow \GG'$ is  a map  of objects $F_0 : G_0 \rightarrow G_0'$ and for any pair of  objects $(A, B)$ in $\GG$, a functor $F_{A,B} : \GG(A,B) \rightarrow \GG'(F_0 A, F_0B)$  preserving composition and identities. That is for all triples of objects $(A, B, C)$  and objects $D$ in $G_0$:
\[
F_{B,C} (-) \circ F_{A,B} (-) = F_{A,C} (- \circ -)
\qquad \text{and} \qquad
iFD(-) = F_{D,D}i (-)
\]
\end{definition}

\begin{definition}
Given two 2-functors $F, E: \GG \rightarrow \GG'$, a \emph{2-natural transformation} $\theta: F \Rightarrow E$ assigns to every object $A$ of $\GG$ a morphism  $\theta A \in G'_1(F_0 A, E_0 A)$ such that for all pairs of objects $(A, B)$ the following generalisation of diagram \eqref{naturaltransformation2} to a diagram in the category of small categories commutes:
\begin{equation}
\label{twonaturaltransformation}
\xymatrix@C=15mm{
\GG(A,B)
	\ar[r]^{F_{A,B}}
	\ar[d]_{E_{A,B}}
& \GG'(F\!A, FB)
	\ar[d]^{\theta B \circ} \\
\GG'(E A, E B)
	\ar[r]_{\circ \theta A}
& \GG'(F\! A, E B)
}
\end{equation}
Here $\circ\theta A$ and $\theta A \circ$ denote pre- and postcomposition with the morphism $\theta A$ on the level of $G_1$ and with $i \theta A$ -- also known as `whiskering with $\theta A$' --  on the level of $G_2$. Thus, for 2-morphisms the condition of the diagram can be written as
\[
\xymatrix@C=15mm{
F\!A	\ar@/^5mm/[r]^{\theta A}_{}="a"
	\ar@/_5mm/[r]_{\theta A}^{}="b"
& EA	\ar@/^5mm/[r]^{Ef}_{}="c"
	\ar@/_5mm/[r]_{Eg}^{}="d"
& EB 
	\ar@{}|{=}[r]
& F\!A	\ar@/^5mm/[r]^{Ff}_{}="e"
	\ar@/_5mm/[r]_{Fg}^{}="f"
& FB	\ar@/^5mm/[r]^{\theta B}_{}="g"
	\ar@/_5mm/[r]_{\theta B}^{}="h"
& EB
\ar@{=>} "a";"b" ^{i\theta A}
\ar@{=>} "c";"d" ^{E\alpha}
\ar@{=>} "g";"h" ^{i\theta B}
\ar@{=>} "e";"f" ^{F\alpha}
}
\]
for morphisms $f, g \in \GG(A,B)$. If  $\theta A$ is an isomorphism for all objects $A$, the 2-natural transformation $\theta$ is called a \emph{2-natural isomorphism}.
\end{definition}

\begin{definition}
2-categories $\GG, \GG'$ are called \emph{2-equivalent} if there are 2-functors
$
F: \GG \longrightarrow \GG'$
and
$E: \GG' \longrightarrow \GG$
such that there are 2-natural isomorphisms 
$\theta: E \circ F \stackrel{\!\sim}{\Longrightarrow} \id_{\GG}$
and
$\sigma: F \circ E \stackrel{\!\sim}{\Longrightarrow} \id_{\GG'}$ .
\end{definition}

%
%
\mysection{Strict 2-groups}
\label{2groups}
\noindent
There are a number of ways to think about a `higher' analogue of a group. A particularly useful and simple one is the  notion of a strict 2-group which can be defined using the language of internal categories.

\begin{definition}
A \emph{strict 2-group} is an internal category in the category of groups.
\end{definition}

This definition of a strict 2-group as a category in the category of groups is equivalent to that of a 2-group as a group in the category of small categories, also known as a categorical group, where the associativity and  unit conditions hold as equalities rather than up to isomorphisms -- hence the `strict' in the name. See 
\cite{math.QA/0307200} for an explanation of both approaches and their equivalence.

\begin{example}
Every group $G$ gives a strict 2-group by taking $G$ as objects and the minimal set of morphisms: $(G, G, \id_G, \id_G, \id_G, \circ)$. 
\end{example}

\begin{example}
Every group $G$ gives a groupoid  with a single object and $G$ as the group of morphisms. If $G$ is abelian, the multiplication  map in $G$ is a morphism and can be used as composition, making the groupoid a strict 2-group.
\end{example}

\begin{example}
For a given group $G$, there is a 2-group with $G_0 = \Aut{G}$, $G_1 = G\ltimes \Aut{G}$, $s(g,F) = F$, $t(g, F) = \Ad{g}\circ F$ -- $\Ad{}$ denoting the adjoint action --, $i(g) = (e,g)$ and $(g', F') \circ (g, F) = (g'g, F)$. This is known as the \emph{automorphism 2-group of $G$}.
\end{example}

The composition being a group morphism is equivalent to the \emph{interchange law} or \emph{middle-four-exchange}, which shows how the multiplication $\cdot$  and composition $\circ$ can be interchanged; for all $(a,b), (a',b') \in  G_1 \!\ _{s}\!\!\times_{t} G_1$:
\begin{equation}
\circ\left(a\cdot a', b\cdot b'\right) 
= \circ\left((a,b) \cdot (a', b')\right)
= \circ(a,b) \, \cdot \,  \circ(a', b')
\label{interchangelaw}
\end{equation}
In particular the interchange law and $sa=tb$ give that\\
\begin{minipage}{85mm}
\begin{equation}
\label{composition}
\begin{aligned}
a\circ b
&= (a\, isa^{-1}\, isa) \circ (isa\, isa^{-1} b) \\
& = (a \circ isa) (isa^{-1} \circ isa^{-1}) (isa \circ b) \quad\\
& = a (isa^{-1}\circ isa^{-1}) b\\
& = a \; isa^{-1} \; b\,. \\
\end{aligned}
\end{equation}
\end{minipage}
\quad
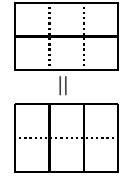
\begin{minipage}{45mm}
\vspace{-3mm}
\begin{displaymath}
\xymatrix@M=0mm@C=4.5mm@R=4.5mm{
\ar@{-}[rrr]\ar@{-}[dd] & \ar@{.}[dd] & \ar@{.}[dd] & \ar@{-}[dd]\\
\ar@{-}[rrr] &&& \\
\ar@{-}[rrr] &\ar@{}[d]^{\;\;\rotateninety{\!\!=}}&& \\
\ar@{-}[rrr] \ar@{-}[dd] & \ar@{-}[dd] & \ar@{-}[dd] & \ar@{-}[dd]\\
\ar@{.}[rrr] &&& \\
\ar@{-}[rrr] &&& \\
}
\end{displaymath}
\end{minipage}\\[2mm]
Here, the second equality is illustrated by the diagram at the right where the inner vertical lines represent multiplication and the middle horizontal line represents composition. The operation of the dotted lines is used first, that of the solid ones second. 
We have thus shown:

\begin{proposition}
To define a strict 2-group we only need to specify the groups $G_0$ and $G_1$ along with the source, target and identity morphisms $s, t, i$. If a composition morphism $\circ$ exists such that the conditions \eqref{categoryconditions} are satisfied, it is uniquely determined by the other maps as $a \circ b = a\,isa^{-1}\, b$.
\end{proposition}

\begin{corollary}
For $a\in \ker s$ we have:
\begin{equation}
\label{compositionismultiplication}
a \circ b
= a \, isa^{-1} \, b 
= a b.
\end{equation}
\end{corollary}

It will prove useful to have the group isomorphism
\begin{equation}
\label{phi}
\phi: G_1 \longrightarrow \ker s \ltimes G_0 
\qquad \qquad 
a \longmapsto (a\, isa^{-1}, sa)
\end{equation}
with inverse $\phi^{-1}(h,g) =  h \, ig$
where the semidirect product comes from the adjoint action  $g\mapsto \Ad{ig}$ of $G_0$ on $\ker s$.
The composition $\circ$ together with $\phi$ lets us move the composition to $\ker s \ltimes G_0$ as $\circ = \phi (\circ) \phi^{-1}\times \phi^{-1}$; for all composable $(k,g), (k', g')\in \ker s \ltimes G_0$:
\begin{align}
\notag
\phi [\phi^{-1}(k,g) \circ \phi^{-1}(k', g')] 
	& = \phi [k \, ig \circ k'ig'] \\
	& = \phi[k\, ig\, is(k\,ig)\, k'ig']
	= \phi[k k' \, ig]
	= (kk', g)
\end{align}
Projecting to the first component with the map $\pi_1$ shows that on the $\ker s$-component composition on $\ker s\ltimes G_0$ is nothing but multiplication -- which is another statement of the corollary above.

We still need morphisms to go along with strict 2-groups. As strict 2-groups are defined as categories in the category of groups, we define their morphisms and 2-morphisms along the same lines to give us a 2-category. 

\begin{definition}
A \emph{morphism of strict 2-groups} is a functor in the category of groups. A \emph{2-morphism of strict 2-groups} is a natural transformation in the category of groups. With these we get the \emph{2-category of strict 2-groups}, $\twoGrp$, as the 2-category of categories, functors and natural transformations in the category of groups.

\end{definition}

%
%
\mysection{Crossed modules}
\label{crossedmodules}
\noindent
Crossed modules appeared in homotopy theory in the mid twentieth century \cite{whitehead-ch2}. The name of the Peiffer identity \eqref{peiffer} used in their definition goes back to the elements defined by Peiffer in her dissertation \cite{peiffer-uizr} which were used by Reidemeister to formulate the equation \cite[(5)]{reidemeister-uir}. 

\begin{definition}
A \emph{crossed module} is a quadruple consisting of two groups, a group homomorphism and an action, $\chi = (G, H, \, \tau \!:\!H\!\rightarrow\!G, \, \alpha\!:\!G \times H \!\rightarrow\!H)$, such that the following dia\-grams commute:
\begin{equation*}
\xymatrix@C=15mm{
	G \times H \ar[d]_{\id \times \tau}  \ar[r]^{\alpha} 
&	H \ar[d]^{\tau} \\
	G \times G \ar[r]_{\Ad{}}
&	G
}
\qquad\qquad\qquad
\xymatrix@C=15mm{
	H \times H \ar[d]_{\tau \times \id}  \ar@/_.0000001pt/[rd]^{\Ad{}} & \\
	G \times H \ar[r]_\alpha 
&	H
}
\end{equation*}
\end{definition}

In other words --  we require $\tau $ to be $G$-equivariant and the action $\alpha$ to satisfy the so-called \emph{Peiffer identity} for all $h,{h'} \in H$:
\begin{equation}
\label{peiffer}
\alpha(\tau h,h') = h{h'}h^{-1}
\end{equation}

\setcounter{example}{0}
\begin{example}
For a group $G$ we get a crossed module $(G, \{e\}, i, \alpha(g, \cdot) = \id_{\{e\}})$ for which  the inclusion $i$ and the identity action are the only possible maps for the given groups.
\end{example}

\begin{example}
For an abelian group $G$ $(\{e\}, G, e, \alpha(e, \cdot) = \id_G)$ is a crossed module where the constant map $e$ and the identity action are the only ones possible.
\end{example}

\begin{example}
For any group $G$  $(\Aut{G}, \,G, \,\tau: g\mapsto \Ad{g}, \,\alpha = \id_{\Aut{G}})$ is a crossed module as $\Ad{F}\,  \tau g = \Ad{F}\, \Ad{g} = \Ad{F(g)} = \tau F(g) = \tau \alpha(\phi,g)$ and $\alpha(\tau g, g') = \Ad{g} (g') = gg'\!g^{-1}$.
\end{example}

\begin{definition}
A \emph{morphism of crossed modules} $(G, H, \tau, \alpha)$, $(G', H', \tau', \alpha')$ is  a pair of group morphisms $(\gamma: G \rightarrow G',\, \delta: H \rightarrow H')$  such that the following diagrams commute:\\
\begin{minipage}{0.5\textwidth}
\begin{equation}
\xymatrix@C=15mm{
	H \ar[d]_{\tau} \ar[r]^\delta
&	H' \ar[d]^{\tau'} \\
	G \ar[r]_{\gamma}
&	G'
}
\stepcounter{equation}
\tag{\theequation a}
\label{crossedmodulemorphism1}
\end{equation}
\end{minipage}
\begin{minipage}{0.5\textwidth}
\begin{equation}
\xymatrix{
	G \times H \ar[d]_<>(0.5){\alpha}  \ar[r]^<>(0.5){\gamma \times \delta}
&	G' \times H' \ar[d]^{\alpha'}\\
	H \ar[r]_{\delta} 
&	H'
}
\tag{\theequation b}
\label{crossedmodulemorphism2}
\end{equation}
\end{minipage}
\end{definition}

\begin{fact}
Crossed modules and morphisms of crossed modules together with component-wise composition and identities give the \emph{category of crossed modules}.
\end{fact}

\begin{definition}
Given morphisms $(\gamma, \delta), (\Gamma, \Delta): (G, H, \tau, \alpha) \rightarrow (G', H', \tau', \alpha')$ of crossed modules, a \emph{2-morphism $\eta: (\gamma , \delta) \Rightarrow (\Gamma, \Delta)$ of crossed modules} is given by the source and target morphisms along with a map $\eta: G\rightarrow H'$ satisfying the condition
\begin{equation}
\label{etacondition}
\eta(g\tilde{g}) = \eta g \, \alpha'(\gamma g , \eta \tilde{g})
\end{equation}
for all $g, \tilde{g}\in G$ as well as for all $g\in G$ and $h\in H$ the conditions\\
\noindent
\begin{minipage}{60mm}
\begin{gather}
\stepcounter{equation}
\tag{\theequation a}
\label{chainhomotopy1}
\tau'  \eta (g) = \Gamma g \: \gamma g^{-1}\\
\tag{\theequation b}
\label{chainhomotopy2}
\eta \tau(h) = \Delta h \: \delta h^{-1} 
\end{gather}
\end{minipage}
\hspace{20mm}
\begin{minipage}{30mm}
\[
\xymatrix@C=20mm{
H 	\ar@<-0.8mm>[r]_\Delta 
	\ar@<0.8mm>[r]^\delta
	\ar[d]_\tau
& H'	\ar[d]^{\tau'}
\\
G 	\ar@<-0.8mm>[r]_\Gamma 
	\ar@<0.8mm>[r]^\gamma
	\ar@/^.0001pt/[ur]^<>(0.25)\eta
& G'
}
\]
\end{minipage}\\
\end{definition}

\begin{remark}
$\eta$ is an ordinary map and \emph{not} required to be a group morphism. Instead, the condition of equation \eqref{etacondition} may look familiar from that of crossed homomorphisms in group cohomology \cite[IV.2]{maclane-h}.  It ensures that the induced map
\[
E: G\longrightarrow H' \ltimes G'
\qquad\qquad
g\longmapsto (\eta g , \gamma g)
\]
is a group morphism. The other two equations are the usual chain homotopy conditions if one considers crossed modules as chain complexes with two non-trivial groups and $\tau$ as their only non-trivial morphism.
\end{remark}

For 2-morphisms  $\eta: (\gamma, \delta) \Rightarrow (\Gamma, \Delta)$, $\bar{\eta} : (\Gamma, \Delta) \Rightarrow (\bar{\Gamma}, \bar{\Delta})$ a (vertical) composition is defined by multiplication in $H'$:
\[
(\bar{\eta} \bullet \eta) : 
(\gamma, \delta) \Longrightarrow (\bar{\Gamma}, \bar{\Delta})
\qquad\qquad
G \longrightarrow H'
\qquad \qquad
g \longmapsto \bar{\eta} g \, \eta g
\]
It is associative as multiplication in $H'$ is and a few computations confirm that $\bar{\eta} \bullet \eta$ is indeed a 2-morphism of crossed modules: For all $g, \tilde{g} \in G$ use the first chain homotopy condition for $\eta$ and the Peiffer identity  to see
\[
\alpha' (\Gamma g, \bar{\eta} \tilde{g}) \, \eta g
	\!\stackrel{\eqref{chainhomotopy1}}{=}\! \alpha' (\tau'\! \eta g \, \gamma g , \bar{\eta} \tilde{g} ) \, \eta g 
	= \alpha' (\tau'\! \eta g, \, \alpha'(\gamma g, \bar{\eta} \tilde{g})) \, \eta g 
	\!\stackrel{\eqref{peiffer}}{=}\! \eta g \, \alpha' (\gamma g, \bar{\eta} \tilde{g}) 
\]
and it follows that
\begin{align*}
(\bar{\eta} \bullet \eta ) \, (g\tilde{g}) 
	= \bar{\eta}( g\tilde{g}) \, \eta(g \tilde{g}) 
&	= \bar{\eta} g \, \alpha'(\Gamma g, \bar{\eta} \tilde{g}) \, \eta g \, \alpha'(\gamma g, \eta \tilde{g}) \\
&	= \bar{\eta} g \, \eta g \, \alpha'(\gamma g, \bar{\eta} \tilde{g}) \,  \alpha'(\gamma g, \eta \tilde{g}) \\
&	= (\bar{\eta} \bullet \eta) g  \, \alpha'(\gamma g, (\bar{\eta} \bullet \eta) \tilde{g}) \,.
\end{align*}
The chain homotopy conditions for $\bar\eta \bullet \eta$ follow straightforwardly:
\begin{gather*}
\tau'  (\bar{\eta} \bullet \eta ) g 
	= \tau' \bar{\eta} g \, \tau' \eta g
	= \bar{\Gamma} g \, \Gamma g^{-1}  \Gamma g \, \gamma g^{-1}
	= \bar{\Gamma} g \, \gamma g^{-1} \\
(\bar{\eta} \bullet \eta) \tau h
	= \bar{\eta} \tau h \, \eta \tau h
	= \bar{\Delta} h \, \Delta h^{-1}  \Delta h \, \delta h^{-1}
	= \bar{\Delta} h \, \delta h^{-1}
\end{gather*}
This, together with the constant maps
\[
i(\gamma, \delta) : 
(\gamma, \delta) \Longrightarrow (\gamma, \delta)
\qquad\qquad
G \longrightarrow H'
\qquad \qquad
g \longmapsto e
\]
as unit morphisms, defines the category of morphisms and 2-morphisms from $\chi$ to $\chi'$, $\XMod (\chi, \chi')$. To get a 2-category we also need  a composition functor which is defined along the same lines as the one given for the horizontal composition of natural transformations in \eqref{horizontalcompositionnaturaltransformations}:
\[
\begin{array}{rccl}
\circ : & \XMod(\chi', \chi'') \times \XMod(\chi, \chi') & \longrightarrow & \XMod(\chi, \chi'')\\
& ((\gamma', \delta') , (\gamma, \delta)) & \longmapsto & (\gamma' \circ \gamma, \delta' \circ \delta) \\
& \!\!\!\!\!(\eta'\!:\! (\gamma', \delta') \!\Rightarrow\! (\Gamma', \Delta'),\;
 \eta \!:\! (\gamma, \delta) \!\Rightarrow \! (\Gamma, \Delta)) & \longmapsto & \Delta' \eta \cdot \eta' \gamma = \eta' \Gamma  \cdot \delta' \eta 
\end{array}
\]
On the object level this is composition in the category of crossed modules, on the morphism level we need to verify that it actually gives a 2-morphism, that is, it satisfies equation \eqref{etacondition}. In a computation similar to that done for the vertical composition we first see that for the horizontal composition with maps named as in the diagram
\[
\xymatrix@C=20mm{
H 	\ar@<0.8mm>[r]^\delta 
	\ar@<-0.8mm>[r]_\Delta
	\ar[d]_\tau
& H'	\ar@<0.8mm>[r]^{\delta'}
	\ar@<-0.8mm>[r]_{\Delta'}
	\ar[d]_{\tau'}
& H''\ar[d]^{\tau''} 		\\
G	\ar@<0.8mm>[r]^\gamma
	\ar@<-0.8mm>[r]_\Gamma	
	\ar@/_.0001mm/[ur]^<>(0.25)\eta
& G'	\ar@<0.8mm>[r]^{\gamma'}
	\ar@<-0.8mm>[r]_{\Gamma'}
	\ar@/_.0001mm/[ur]^<>(0.25){\eta'}
& G''
}
\]
we have for all $g, \tilde{g} \in G$
\begin{gather}
\notag
\Delta' \alpha'(\gamma g, \eta \tilde{g}) \, \eta'\! \gamma g
\stackrel{\eqref{crossedmodulemorphism2}}{=}
	 \alpha'' (\Gamma' \! \gamma g , \Delta' \! \eta \tilde{g} ) \, \eta' \! \gamma g 
\stackrel{\eqref{chainhomotopy1}} { =}
 	\alpha'' (\tau'' \! \eta' \! \gamma g \: \gamma' \! \gamma g, \Delta' \! \eta \tilde{g}) \,\eta' \! \gamma g \\
 =	\alpha'' (\tau'' \! \eta' \! \gamma g, \alpha'' (\gamma' \! \gamma g, \Delta' \! \eta \tilde{g})) \, \eta'\! \gamma g 
\stackrel{\eqref{peiffer}}{=}
	 \eta' \! \gamma g \: \alpha'' (\gamma' \! \gamma g, \Delta' \! \eta \tilde{g})
\label{xmodhorizontalcompositionaux}
\end{gather}
With this it follows that
\begin{align*}
(\eta'\! \circ \eta) (g\tilde{g})
&	= \Delta' \eta (g\tilde{g}) \: \eta'\! \gamma (g\tilde{g}) \\
&	= \Delta' \eta g \;\; \Delta' \alpha' (\gamma g, \eta \tilde{g}) \; \eta'\!\gamma g \;\; \alpha''(\gamma'\! \gamma g, \eta'\! \gamma g) \\
\formularemark{\eqref{xmodhorizontalcompositionaux}} &	= \Delta' \eta g \; \eta'\! \gamma g \; \alpha'' (\gamma'\!\gamma g, \Delta' \eta \tilde{g}) \: \alpha''(\gamma'\! \gamma g, \eta'\! \gamma g) \\
&	= \Delta' \eta g \; \eta'\! \gamma g \; \alpha'' (\gamma'\!\gamma g, \Delta' \eta \tilde{g} \; \eta'\! \gamma \tilde{g}) \\
&	= (\eta'\! \circ \eta) g  \; \alpha''(\gamma' \! \gamma g, (\eta' \! \circ\eta) \tilde{g}) \,.
\end{align*}	
Again, the chain homotopy conditions are verified easily,
\begin{gather*}
\tau''(\eta' \! \circ\eta) \, g 
	= \tau'' \! \Delta' \eta g \; \tau''\!\eta'\!\gamma g 
	\stackrel{\textrm{(\ref{crossedmodulemorphism1},\ref{chainhomotopy1})}}{=} 
	\Gamma'\!\tau'\! \eta g \; \Gamma'\!\gamma g \: \gamma' \!\gamma g^{-1}
	\stackrel{\eqref{chainhomotopy1}}{=} \Gamma' \Gamma g \; \gamma'\! \gamma g^{-1}\\
(\eta'\! \circ \eta) \, \tau h 
	= \Delta' \eta \tau h \; \eta' \! \gamma \tau h
	\stackrel{\textrm{(\ref{chainhomotopy2}, \ref{crossedmodulemorphism1})}}{=} \Delta' \! \Delta h \: \Delta'\! \delta h^{-1} \, \eta'\! \tau'\! \delta h
	\stackrel{\eqref{chainhomotopy2}}{=} \Delta' \! \Delta h \; \delta'\! \delta h^{-1}\,,
\end{gather*}
showing that the definition of horizontal composition above gives another 2-morphism. This composition is associative as the multiplication in $H''$ is: 
\begin{multline*}
\eta''\! \circ (\eta'\! \circ \eta)
	= \eta ' \! \circ ( \Delta' \eta \cdot \eta' \! \gamma)
	= \Delta'' (\Delta' \eta \cdot \eta' \! \gamma) \cdot \eta'' \! \gamma'\! \gamma\\
	= \Delta''\! \Delta' \eta \cdot \Delta'' \eta' \! \gamma \cdot \eta'' \! \gamma' \! \gamma 
	= \Delta''\! \Delta' \eta \cdot (\eta''\! \circ \eta') \gamma
	= (\eta''\! \circ \eta') \circ \eta
\end{multline*}
Finally, the unit functor for a crossed module $\chi$ is defined by
\[
i\chi: 1 \longrightarrow \XMod(\chi, \chi) 
\qquad \quad
* \longmapsto (\id_{G}, \id_{H})
\qquad \quad
\id_{*} \longmapsto i(\id_{G}, \id_{H})
\]
On the object level this is just the unit of the category of crossed modules and crossed module morphisms, on the morphism level, we easily verify for any 2-morphism $\eta:(\gamma, \delta) \Rightarrow (\Gamma, \Delta): \chi \rightarrow \chi'$ that the unit functor satisfies:
\[
\eta\circ i\chi(\id_*)
	= \Delta i\chi(\id_*) \; \eta \id_G 
	= \eta
\qquad
i\chi'(\id_*) \eta
	= \id_{H'} \eta \; i\chi'(\id_*) \gamma = \eta
\]
Hence we have shown:
\begin{propositiondefinition} 
Crossed modules, crossed module morphisms and crossed module 2-morphisms together with the structures defined above form the \emph{2-category $\XMod$ of crossed modules}.
\end{propositiondefinition}

%
%
\mysection{Strict 2-groups give crossed modules}
\label{2groupstomodules}
\noindent
We want to define a 2-functor $T$ from the 2-category $\twoGrp$ to the 2-category $\XMod$. 
For this we need a map $T_0$ that gives us a crossed module for each strict 2-group and a functor $T_{\GG, \GG'}$ for any two 2-groups $\GG$, $\GG'$.

On the \emph{object level} define
\[
T_0(\GG)
	= (G_0,\, 
		\ker s, \,
		\tau = t\restrictedto{\ker s}, \,
		\alpha(g,h) = \Ad{ig} h
	)
\]	
To see that this is a crossed module, note that all components are groups and group morphisms as required. Equivariance of $\tau$ follows from $t$ being a group homomorphism with $t\circ i = \id_G$:
\[
\tau \alpha(g,h)
= t(ig \, h \, ig^{-1})
= g\,th\,g^{-1}
= \Ad{g}\,\tau (h)
\]

The Peiffer identity only concerns elements of $H=\ker s$. Their role in $G_1$ is more visible in the isomorphic group $\ker s\ltimes G_0$ after using the isomorphism $\phi$ introduced in \eqref{phi}. Given $a\in G_1$, write $\ti{a}$ for the first component of $\phi(a)$:
\[
\pi_1 \phi(a) 
	= \pi_1(a \, isa^{-1}, sa)
	= a\, isa^{-1}
	=: \ti{a}
\]
With this the interchange law  gives:
\[
\hspace{-2mm}
\begin{array}{ccc}
\pi_1\!\left[[(\ti{a}, sa) (\ti{a}', sa')] \! \circ \! [(\ti{b}, sb) (\ti{b}', sb')]\right]
\hspace{-4mm}
& \stackrel{\eqref{interchangelaw}}{=}& 
\hspace{-3mm}
\pi_1\!\left[[(\ti{a}, sa) \!\circ\! (\ti{b}, sb)] [(\ti{a}', sa') \!\circ\! (\ti{b}',sb')]\right] 
\\[2mm]
\rotateninety{=}&& \rotateninety{=}
\\
\pi_1\left[(\ti{a}\, isa \, \ti{a}'\, isa^{-1}, sa\,sa') 
		\circ (\ti{b}\, isb \, \ti{b}' \, isb', sb\, sb')\right] 
\hspace{-21mm}
&&
\qquad\pi_1\left[(\ti{a}\ti{b}, sb)(\ti{a}'\ti{b}', sb')\right]
\\
\rotateninety{=}&& \rotateninety{=}
\\
\pi_1\left[\ti{a}\, isa \, \ti{a}' \, isa^{-1} \, \ti{b} \, isb\, \ti{b}' \, isb^{-1}, sb\,sb'\right]
\hspace{-5mm}
&&
\pi_1\left[\ti{a}\ti{b}\, isb\, \ti{a}' \ti{b}' \, isb^{-1}, sb\, sb'\right]
\\
\rotateninety{=}&& \rotateninety{=}
\\
\ti{a}\,\alpha(sa, \ti{a}')\,\, \ti{b} \, \alpha(sb, \ti{b}')
&&
\ti{a}\ti{b}\, \alpha(sb, \ti{a}'\ti{b}')
\end{array}
\]
Cancelling $\ti{a}$ on the left, $\ti{b}\, \alpha(sb , \ti{b}')$ on the right and using that $sa = tb$ gives
\[
\alpha(tb,\ti{a}') = \ti{b}\, \alpha(sb, \ti{a}')\, \ti{b}^{-1},
\]
which, for $h=\ti{b}=b\in \ker s$ and $h' = \ti{a}' = a'\in \ker s$, yields the Peiffer identity \eqref{peiffer}: $\alpha(\tau h,h') = hh'h^{-1}$. It follows that $(G_0, \ker s, t\restrictedto{\ker s}, \alpha)$ is a crossed module as claimed.

Next define the functors
\[
T_{\GG, \GG'} = (T_1, T_2) : \twoGrp(\GG, \GG') \longrightarrow \XMod(T_0 \GG, T_0 \GG')\,.
\]
On the object level, that is  \emph{a map $T_1$ of strict 2-group morphisms to crossed module morphisms}
\[
T_1: (f_0, f_1) \longmapsto (\gamma = f_0, \delta = {f_1}\restrictedto{\ker s})
\]
which is well-defined as functoriality of $(f_0, f_1)$ implies $f_1(\ker s) \subset \ker s'$.
With $(f_0, f_1)$ being a functor in the category of groups we have  for all $g\in G$ and $h\in H$ that
\begin{gather*}
\gamma(\tau h) 
	= f_0 (th) 
	= t'({f_1}\restrictedto{\ker s}(h)) 
	= \tau'  \delta h \, , \\
\delta(\alpha(g,h)) 
	= {f_1}\restrictedto{\ker s}(ig \,  h \, ig^{-1})
	= f_1(ig) \, f_1(h) \,  f_1(ig^{-1}) = \qquad\\
	\qquad \qquad
	= if_0(g) \, f_1(h) \, if_0(g)^{-1} 
	= \alpha'(f_0(g), {f_1}\restrictedto{\ker s}(h))
	= \alpha'(\gamma g,\delta h),
\end{gather*}
proving that $T_1$ indeed gives a morphism of crossed modules.

The two maps $T_0$ and $T_1$ we defined so far define a functor as two further computations verify:
\[
T_1 i \GG
	= T_1 {(\id_{G_0}, \id_{G_1})}
	= (\id_{G_0}, {\id_{G_1}}\restrictedto{\ker s})
	= i\, T_0 \GG
\]
\begin{align*}
T_1((f_0', f_1') \circ  (f_0, &\, f_1)) 
	= T_1(f_0'\circ f_0, \, f_1' \circ f_1)
	= (f_0' \circ f_0, \, f_1' \circ {f_1}\restrictedto{\ker s}) =  \\
&	= (f_0', {f_1'}\restrictedto{\ker s'}) \circ (f_0, {f_1}\restrictedto{\ker s}) 
	= T_1(f_0', f_1') \circ T_1(f_0, f_1)
\end{align*}

Now consider the morphism component of the functor $T_{\GG, \GG'}$ --  a map $T_2$ on the \emph{level of 2-morphisms}. For two 2-group morphisms $F=(f_0, f_1)$, $E=(e_0, e_1): \GG \rightarrow \GG'$ of strict 2-groups and a 2-morphism $F\Rightarrow E$ between those morphisms, that is a natural transformation which is given by a group morphism
$\theta: G_0 \rightarrow G_1'$,
we define a 2-morphism $T_1 F \Rightarrow T_1 E$ of crossed modules. To do that we need a map $\eta: G=G_0 \rightarrow H'$. Composing $\theta$ with the isomorphism $\phi'$ and the projection $\pi_1'$ gives
\[
\eta: G = G_0 \stackrel{\theta}{\longrightarrow}
G_1' \stackrel{\phi'}{\longrightarrow}
\ker s' \ltimes G_0' \stackrel{\pi_1'}{\longrightarrow}
\ker s' = H' \,.
\]
With this we define
\[
T_2:  [\,\theta: F \Rightarrow E \quad G_0 \rightarrow G_1'\,]\
\longmapsto 
[\,\eta: T_1 F \Rightarrow T_1E \quad \eta = \pi_1 \phi' \theta: G \rightarrow H'\,]\,.
\]
In particular we have for all $g\in G$ that
\[
\eta g 
= \pi_1' \phi' \theta g
= \pi_1' ( \theta g \: i'\!s'\!\theta g^{-1}\!, s'\!\theta g)
= \theta g \: i'\! s'\! \theta g^{-1}.
\]
For $g, \tilde{g} \in G$ using the previous computation, the semidirect product, and $s'\!\theta = \gamma$ gives
\begin{align*}
\eta(g\tilde{g})
&	 = \pi_1' ( \phi' \theta g \cdot \phi' \theta \tilde{g} 
) \\
& 	= \pi_1' [(\theta g \, i'\!s'\! \theta g^{-1}, s'\!\theta g) ( \theta \tilde{g}\, i'\! s'\! \theta \tilde{g}^{-1}, s' \! \theta \tilde{g})] \\
& 	= \theta g \, i'\! s'\! \theta g^{-1} \, \alpha'(s'\!\theta g, \theta \tilde{g} \, i'\!s'\! \theta \tilde{g}^{-1}) \\
& = \eta g \, \alpha'(\gamma g, \eta \tilde{g})
\end{align*}
thus proving that the condition from equation \eqref{etacondition} in the definition of a 2-morphism of crossed modules holds. The first chain homotopy property \eqref{chainhomotopy1} follows from applying the natural transformation $\theta$ to $g\in G$, giving morphisms $\theta(g): f_0(g) \rightarrow e_0(g)$:
\begin{align*}
\tau'  \eta g 
	& = t'\restrictedto{\ker s'} \pi_1' \phi' \theta g 
	 = t' (\theta g \: i'\! s'\! \theta g^{-1}) 
	 = t'\theta g \: t'\!i'\! s'\! \theta g^{-1} \\
	& = t'\! \theta g \, s'\! \theta g^{-1} 
	 = e_0 g \, f_0  g^{-1} 
	 = \Gamma g \, \gamma g^{-1}	
\end{align*}
\begin{minipage}{60mm}
The first chain homotopy property \eqref{chainhomotopy2} follows from using the natural transformation diagram of $\theta$ for a morphism
\[
(h: e \rightarrow \tau h = th) \in H,
\]
seen on  the right. We conclude that
\end{minipage}
\qquad
\begin{minipage}{50mm}
\[
\vspace{3mm}
\xymatrix@C=20mm{
\hspace{-14mm}f_0(e) = e 	
	\ar[r]^<>(.5){f_1(h) = \delta h} 
	\ar[d]_{ie = \theta e}
& f_0(t h) 
	\ar[d]^{\theta th} 
\\
\hspace{-14mm}e_0(e) = e
	\ar[r]_<>(.5){e_1(h) = \Delta h}
& e_0(t h)
}
\]
\end{minipage}\\
\begin{align*}
\eta \tau (h)
	& = \pi_1 \phi' \theta th
	= \theta t h \: i'\!s'\! \theta t h^{-1} \\
	\formularemark{$th^{-1}= e$,  compose with identity}
	& = (\theta t h \, ie) \circ (\delta h \, \delta h ^{-1}) \\
	\formularemark{interchange law \eqref{interchangelaw}}
	& = (\theta t h \circ \delta h) (ie \circ \delta h^{-1}) \\
	\formularemark{natural transformation diagram}
	& = (\Delta h \circ i_e) \, \delta h^{-1}
	= \Delta h \, \delta h^{-1}\,,
\end{align*}
completing the proof that $T_2$ is a map from 2-morphisms of strict 2-groups to 2-morphisms of crossed modules.

We now verify that $T_{\GG, \GG'} = (T_1, T_2)$ is a functor as claimed. It preserves identities for any $F: \GG \rightarrow \GG'$:
\[
i \, T_1 F 
	\stackrel{\eqref{identitynaturaltransformation}}= [\,T_1 F \Rightarrow T_1 F \quad g \mapsto e\,]
	= T_2 (F \Rightarrow F \quad g \mapsto i F g)
	= T_2 \, iF
\]
as well as vertical composition; where we have for vertically composable 2-morphisms $\bar\theta, \theta: \GG \rightarrow \GG'$:
\begin{align*}
T_2 (\bar\theta \bullet \theta : g \mapsto (\bar \theta \bullet \theta ) g )
& 	= [\, g \mapsto \pi_1 \phi' (\bar \theta g \circ \theta g) ]\\
	\formularemark{\eqref{compositionismultiplication}} 	
&	= [\, g \mapsto \pi_1 \phi' (\bar\theta g \, is\theta g ^{-1} \, \theta g)]\\
&	= [\, g \mapsto \bar \theta g \, is\bar\theta g^{-1} \, \theta g \; is g^{-1} \, isis\bar\theta g \, is \bar\theta g^{-1}] \\	
&	= T_2(\bar \theta) \bullet T_2(\theta)\,. 
\end{align*}
This completes the proof that $T: \twoGrp \rightarrow \XMod$ is a 2-functor.

%
%
\mysection{Crossed modules give strict 2-groups}
\label{modulesto2groups}
\noindent 
In this section we define a 2-functor $S: \XMod \rightarrow\twoGrp$. 
On the \emph{object level} we begin with  a crossed module $\chi = (G, H, \, \tau \!:\!H\!\rightarrow\!G, \, \alpha\!:\!G \times H \rightarrow H)$ and define groups 
\[
G_0 = G 
\qquad \textrm{and} \qquad
G_1 = H \ltimes G
\]
with the semidirect product's structure given by the action $\alpha$. Furthermore we define the source, target and identity maps
\begin{gather*}
s: H \ltimes G \longrightarrow G
\quad
(h,g) \longmapsto g
\qquad\qquad
t: H \ltimes G \longrightarrow G
\quad
(h,g) \longmapsto \tau h \, g \\
i: G \longrightarrow H \ltimes G 
\quad
g\longmapsto (e,g) .
\end{gather*}
$s$ is a projection of the semidirect product to its second factor $G$ and $i$ is the inclusion of $G$ in the semidirect product. Hence both are morphisms. The fact that $t$ is a morphism follows from equivariance  of $\tau$:
\begin{align*}
t((h,g)\cdot(h', g') ) 
&= t(h\, \alpha(g,h'), gg')
= \tau h \, \tau (\alpha(g,h'))\,gg' \\
& = \tau h \, g \, \tau h' g^{-1} g g'
= t(h,g) \cdot t(h', g')
\end{align*}
The proposition in section~\ref{2groups} asserts that the composition of a strict 2-group with the maps we already have has to be defined as
\[
(j,\tau h \, g) \circ (h, g)
	= (j, \tau h\, g) (e, \tau h \, g)^{-1} (h,g)
	= (j, e) (h,g )
	= (jh, g)\,.
\]
We need to check that this is a morphism, that is, the interchange law holds:
\begin{align*}
\circ \, [(j, \tau h\, g) , \,& (h,g) ]  \cdot \circ [ (j', \tau h' g') , \, (h', g') ] \\
& = (jh,g)  \cdot (j'h', g') 
= (jh \,\alpha(g,j'h'), \, gg') \\
& = (jh \,\alpha(g,j')\, h^{-1} h \,\alpha(g,h'),\, gg') \\
	\formularemark{Peiffer \eqref{peiffer}}
& = (j \,\alpha(\tau h, \alpha(g,j'))\,  h \, \alpha(g,h'),\, gg') \\
& = ( j \,\alpha(\tau h \, g,j')\, h \,\alpha(g,h'), gg' ) \\
& = \circ [ (j \,\alpha(\tau h\, g,j'),\, \tau h\, g \, \tau h' g' ), (h \,\alpha(g,h'), gg' )] \\
& = \circ [ (j, \tau h \, g)  \cdot (j', \tau h' g')\, , \, (h,g) \cdot (h', g') ]\\
\end{align*}
A few simple computations show that the remaining conditions for a strict 2-group are satisfied as well:
\begin{gather*}
s\circ i(g) = s(e,g) = g
\qquad 
t\circ i(g) = t(e,g) = \tau e\, g = g\\
(h,g) \circ i(s(h,g)) = (h,g) \circ (e,g) = (h,g)\\
i(t(h,g)) \circ (h,g) = (e, \tau h \, g) \circ (h,g) = (h,g)\\
(k, \tau j \, \tau h \, g) \circ ((j, \tau h \, g) \circ (h,g)) = (k, \tau j \, \tau h \, g) \circ (jh,g) = (k(jh), g) = \qquad \\
\quad\qquad = ((kj)h, g) = (kj, \tau h \, g) \circ (h,g) = ((k,\tau j \, \tau h \, g) \circ (j, \tau h \, g)) \circ(h,g)\\
s((j, \tau h \, g) \circ (h,g)) = s(jh, g) = g = s(h,g)\\
t((j, \tau h \, g) \circ (h,g)) = t(jh ,g ) = \tau j \, \tau h \, g = t(j, \tau h \, g)
\end{gather*}
With that being established we can assign a strict 2-group to each crossed module, defining $S$ on objects of $\XMod$:
\[
S_0: (G,H, \tau, \alpha) \longmapsto (G, H\ltimes G, s, t, i)
\]

For crossed modules $\chi$, $\chi'$ we now define the functor 
\[
S(\chi, \chi') : \XMod(\chi, \chi') 
\longrightarrow 
\twoGrp(S_0(\chi), S_0(\chi'))\,.
\]
On objects of $\XMod(\chi, \chi')$, that is on \emph{morphisms} of crossed modules, we define the map
\[
S_1 : [\,(\gamma, \delta) : \chi \rightarrow \chi'\,]
\longmapsto 
(f_0 = \gamma, \, f_1 = (\delta, \gamma))\,.
\]
It gives a pair of group morphisms and we have to verify that this pair is a morphism in  \twoGrp, that is, an internal functor in the category of groups. We see that it preserves identity maps for all $g\in G_0$
\[
f_1 (ig) 
	= (\delta , \gamma) (e,g)
	= (e, \gamma g)
	= i f_0(g)
\]
as well as respecting composition of morphisms $(h,g), (h', g') \in  \ker s \ltimes G_0$ 
\begin{gather*}
(\delta, \gamma) ((h',g') \circ (h,g))
	= (\delta,\gamma) (h'h, g)
	= (\delta h' \delta h, \gamma g) \qquad \\
\qquad \qquad
	= (\delta h' , \gamma g') \circ (\delta h , \gamma g)
	= (\delta, \gamma)(h', g') \circ  (\delta, \gamma)(h,g).
\end{gather*}

Next we  verify that the  $S_1$ preserve identities
\[
S_1 i{(G,H, \tau, \alpha)}
	= T_1 (\id_G, \id_H)
	= (\id_G, \id_{H \ltimes G})
	= i{S_0 (G, H, \tau, \alpha)}
\]
and are compatible with the compositions of morphisms in $\XMod$ and $\twoGrp$, 
\begin{gather*}
S_1 ((\gamma', \delta') \circ (\gamma,\delta))
	= S_1(\gamma'\!\gamma, \delta'\!\delta)
	= (\gamma'\!\gamma, (\delta'\!\delta, \gamma'\!\gamma)) \qquad\\
\qquad\qquad
	= (\gamma', (\delta', \gamma')) \circ (\gamma, (\delta, \gamma))
	= S_1(\gamma', \delta') \circ S_1(\gamma,\delta)\,.
\end{gather*}

On the level of morphisms in $\XMod(\chi, \chi')$, that is \emph{2-morphisms} between morphisms $(\gamma, \delta), (\Gamma, \Delta): \chi \rightarrow \chi'$ of crossed modules, we define
\[
S_2: [\, \eta: (\gamma, \delta) \Rightarrow (\Gamma, \Delta) \quad G\rightarrow H'\,] 
	\longmapsto
	[\, G \rightarrow H'\ltimes G' \quad g\mapsto (\eta g, \gamma g) ]
\]
\noindent
\begin{minipage}{70mm}
To verify that this gives a 2-morphism of 2-groups, that is an internal natural transformation in the category of groups $S_1(\gamma, \delta) \Rightarrow S_1(\Gamma, \Delta)$, we need to show that the diagram at the right commutes for all 2-group morphisms $(h,g): g \rightarrow \bar g$.\\[-3mm]
\end{minipage}
\quad
\begin{minipage}[c]{50mm}
\[
\vspace{2mm}
\xymatrix@C=20mm{
\gamma g	\ar[r]^{(\delta h, \gamma g)}
		\ar[d]_{S_2 \eta (g)}	
& \gamma \bar g
		\ar[d]^{S_2 \eta (\bar g)}
\\
\Gamma g \ar[r]_{(\Delta h, \Gamma g)}
& \Gamma \bar g
}
\]
\end{minipage}

\noindent
To show this we first note that the Peiffer identity gives:
\begin{align*}
\eta \bar g \, \delta h
&	= \eta (\tau h\, g) \,  \delta h 
	= \eta \tau h \: \alpha'( \gamma \tau h, \eta g) \, \delta h  \\
\formularemark{Peiffer \eqref{peiffer}}
&	= \Delta h \: \alpha'(\tau' \delta h^{-1}\!,\, \alpha'(\gamma \tau h, \eta g))
	= \Delta h \: \alpha'(\tau' \delta h^{-1} \gamma\tau h, \eta g)\\
\formularemark{\eqref{crossedmodulemorphism1}}	
&	= \Delta h \: \alpha'(\gamma \tau h^{-1} \gamma \tau h, \eta g)
	= \Delta h \: \alpha'(e, \eta g ) 
	= \Delta h \: \eta g \,,
\end{align*}
which implies the fact we want to show by re-writing the composition in $S_0 \chi'$ in terms of the multiplication \eqref{composition}:
\begin{align*}
S_2 \eta \bar g \circ (\delta h, \gamma g)
&	= (\eta \bar g, \gamma \bar g) \circ (\delta h , \gamma g) 
	= (\eta \bar g, \gamma \bar g ) (e, \gamma \bar g^{-1}) (\delta h, \gamma g) \\
&	= (\eta \bar g \, \delta h, \gamma g)
	= (\Delta h \, \eta g, \gamma g)
	= (\Delta h, \Gamma g) (e, \Gamma g^{-1}) (\eta g, \gamma g) \\
&	= (\Delta h, \Gamma g) \circ ( \eta g, \gamma g) 
	= \Delta h \circ S_2 \eta g \,.
\end{align*}

\noindent
Furthermore $S_2$ preserves units
\[
S_2 [\,i{(\gamma, \delta)}: g\mapsto e\,] 
	= [\,\theta: g\mapsto (e,\gamma g) = i' \gamma g = i{S_1(\gamma, \delta)}(g)] 
	= i{S_1 (\gamma, \delta)} \,,
\]
vertical compositions thanks to composition of morphisms in the 2-group being multiplication on $\ker s'$ \eqref{compositionismultiplication}
\begin{align*}
S_2(\bar \eta \bullet \eta) 
	= S_2(g \mapsto \bar \eta g \, \eta g) 
&	= [\, g \mapsto (\bar \eta g \, \eta g, \gamma g) ] \\
& 	= [\, g \mapsto (\bar \eta g, \Gamma g) \,] 
			\bullet [\, g \mapsto (\eta g, \gamma g)  ] \\
&	= S_2(\bar \eta) \bullet S_2(\eta) \,,
\end{align*}
and horizontal compositions pretty much by using definitions
\begin{align*}
S_2( \eta' \circ \eta ) 
	 = S_2(\eta' \gamma \cdot \delta' \eta)
	& = [\, g \mapsto (\eta' \Gamma g \cdot \delta'\! \eta g, \gamma'\!\gamma g)] \\
	& = [\, g' \mapsto (\eta'\! g', \gamma'\! g) ] \circ [ \, g\mapsto (\eta g , \gamma g) ] \\
	& = S_2 (\eta') \circ S_2(\eta)
\end{align*}
-- completing the proof that  $S = (S_0, S_1, S_2)$ is a 2-functor.

%
%
\mysection{Round Trip}
\label{roundtrip}
\noindent
With the 2-functors $S$ and $T$ defined in the previous sections it is reasonable to hope that the 2-categories of strict 2-groups and crossed modules are 2-equivalent. To show this we need to find 2-natural isomorphisms
\[
\xi : ST \stackrel{\sim}{\Longrightarrow} \id_{\twoGrp}
\qquad \qquad \text{and} \qquad \qquad
\zeta: TS \stackrel{\sim}{\Longrightarrow} \id_{\XMod} \,.
\]
We begin by spelling out what the 2-functor $ST: \twoGrp \rightarrow \twoGrp$ does:
\begin{gather*}
(G_0, G_1, s, t, i)
	\stackrel{T_0}{\longmapsto}
	(G_0, \ker s, t\restrictedto{\ker s}, \Ad{i})
	\stackrel{S_0}{\longmapsto}
	(G_0, \ker s \ltimes G_0, \pi_2, t\pi_1 \!\cdot\! \pi_2, i_2) 
\\
[\, F = (f_0, f_1) : \GG \rightarrow \GG' \,] 
	\stackrel{T_1}{\longmapsto}
	(f_0, {f_1}\restrictedto{\ker \pi_2})
	\stackrel{S_1}{\longmapsto}
	(f_0, ({f_1}\restrictedto{\ker \pi_2}, f_0))
\\
[\, \theta : F \Rightarrow E \,]
	\stackrel{T_2}{\longmapsto}
	[\, g \mapsto \pi_1' \phi' \theta g\,]
	\stackrel{S_2}{\longmapsto}
	[\, g \mapsto (\pi_1' \phi' \theta g, f_0 g) ]
\end{gather*}
For any strict 2-group $\GG$ define the 2-group isomorphism $\xi \GG$ by
\[
\xi \GG = (\id_{G_0}, \phi^{-1}) : TS \GG \longrightarrow \GG
\]
where $\phi^{-1}(h,g) = h \, ig$ is the inverse of the isomorphism \eqref{phi}. To confirm this gives a natural transformation we need to check commutativity of diagram \eqref{twonaturaltransformation} which for our 2-functors looks like this:
\[
\xymatrix@C=15mm{
\twoGrp (\GG, \GG') 
	\ar[r]^<>(.5){(ST)_{\GG, \GG'}}
	\ar[d]_{\id_{\GG, \GG'}}
&\twoGrp(ST\GG, ST\GG')
	\ar[d]^{\xi \GG' \circ}
\\
\twoGrp(\GG, \GG')
	\ar[r]_<>(.5){\circ \xi \GG}
&\twoGrp(ST \GG, \GG')
}
\]
We then verify that for objects $(f_0, f_1) : \GG \rightarrow \GG'$ the diagram commutes by using that any morphism $h$ in the 2-group can be uniquely written as $k\, isg$ where $k\in \ker s$ and $g = sh$:
\begin{align*}
\xi \GG' \circ (ST)_{\GG, \GG'} (f_0, f_1) (g, k\, isg)
& 	= (\id_{G_0}, \phi'^{-1}) \circ (f_0 g, ({f_1}\restrictedto{\ker \pi_2} k, f_0 g)) \\
&	= (f_0 g, f_1k \, i f_0 g)  = (f_0 g, f_1k\, f_1 ig) \\
&  	= (f_0 g, f_1 (k\, ig)) 
	= (f_0, f_1) \circ (g, k \, ig) \\
&	= \id_{\GG, \GG'}(f_0, f_1) \circ (\id_{G_0}, \phi^{-1}) (g, (k,g))
\end{align*}
The verification for 2-morphisms $\theta: F\Rightarrow E \quad \GG\rightarrow \GG'$ is done by expressing composition through multiplication and keeping in mind -- as well as carefully distinguishing -- the definitions of the various compositions:
\begin{align*}
i \xi \GG' \circ (ST)_{\GG, \GG'} (\theta) (g)
&	= \phi^{-1} (\pi_1' \phi' \theta g , f_0 g) \\
&	= \theta g\, is\theta g^{-1}  if_0 g 
	= \theta g\, is \theta g^{-1}  f_1 i g \\
\formularemark{\eqref{composition}} 	
&	= \theta g \circ f_1 i g = \theta \xi \GG  g \circ f_1 i g \\
\formularemark{\eqref{horizontalcompositionnaturaltransformations}} 
& 	= (\theta \circ  i \xi \GG) (g) 
	= (\id_{\GG, \GG'}(\theta) \circ i\xi \GG) (g)	
\end{align*}

The other direction is  more straightforward. Again, we make explicit the 2-functor $TS: \XMod \rightarrow \XMod$:
\begin{gather*}
(G, H, \tau, \alpha) 
	\stackrel{S_0}{\longmapsto}
	(G, H \ltimes G, \pi_2, \tau \pi_1 \cdot \pi_2, i_2)
	\stackrel{T_0}{\longmapsto}
	(G, \ker \pi_2, \tau {\pi_1}\restrictedto{\ker \pi_2}, \Ad{i_2})
\\
(\gamma, \delta) 
	\stackrel{S_1}{\longmapsto}
	(\gamma, (\delta, \gamma))
	\stackrel{T_1}{\longmapsto}
	(\gamma, (\delta, \gamma)\restrictedto{\ker \pi_2})
\\
[\eta: (\gamma, \delta) \Rightarrow (\Gamma, \Delta)]
	\stackrel{S_2}{\longmapsto}
	[\, g \mapsto (\eta g, \gamma g) ]
	\stackrel{T_2}{\longmapsto}
	[\, g \mapsto (\pi_1' \phi' (\eta g, \gamma g), \gamma g) = (\eta g, e)]
\end{gather*}
For any crossed module $\chi$ we define the crossed module isomorphism 
\[
\zeta \chi = (\id_G, {\pi_1}\restrictedto{\ker \pi_2}) : TS \chi \longrightarrow \chi
\]
and verify that the diagram of natural transformations \eqref{twonaturaltransformation} commutes: On the level of objects in $\XMod(\chi, \chi')$ we get for a crossed module morphism $(\gamma, \delta): \chi \rightarrow \chi'$ that 
\[
\zeta \chi' \circ TS (\gamma, \delta)
	= (\gamma, \delta {\pi_1}\restrictedto{\ker \pi_2})
	= (\gamma, \delta) \circ (\id_G, {\pi_1}\restrictedto{\ker \pi_2})
	= \id_{\chi, \chi'} (\gamma, \delta) \circ \zeta \chi
\]
On the morphism level, for crossed module 2-morphisms $\eta: (\gamma, \delta) \Rightarrow (\Gamma, \Delta)$ we have
\begin{multline*}
i \zeta \chi' \circ TS (\eta)
	= i \zeta \chi' \circ [\, g \mapsto (\eta g, e) ] 
	= \pi_1 [\, g \mapsto (\eta g, e) ]
	= \eta \\
	= \Delta e \cdot \eta \id_G
	= \eta \circ i\zeta \chi
	= \id_{\chi, \chi'} (\eta)  \circ i \zeta \chi \,.
\end{multline*}
\noindent
We have thus shown:
\begin{theorem}
The 2-categories $\twoGrp$ and $\XMod$ are 2-equivalent.
\end{theorem}

This equivalence enables us to get the best of both worlds and go back and forth between strict 2-groups and crossed modules when working with them. Strict 2-groups are well known objects in category theory and lend themselves for category theoretical reasoning, while crossed modules are frequently more convenient for computations and appear in examples from homotopy theory or when working with central extensions as is detailed in the following section.

%
%
\mysection{Central extensions give crossed modules}
\noindent
A class of examples for crossed modules and thus strict 2-groups can be constructed from central extensions.
\begin{definition}
For a group $G$, a \emph{central extension} is a group homomorphism $\tau: H\rightarrow G$ such that $\ker(\tau)$ is in the centre of $H$ and
\[
1 \longrightarrow \ker(\tau) \longrightarrow H \stackrel{\tau}{\longrightarrow} G \longrightarrow 1
\]
is a short exact sequence.
\end{definition}

Given a central extension $\tau:H\rightarrow G$ all that is missing to the quadruple we need for a crossed module is an action of $G$ on $H$.  To construct that, choose an arbitrary section $s: G\rightarrow H$ of $\tau$ and define a morphism 
\[
\alpha: G \times H \longrightarrow H
\qquad \qquad
(g, h )  \longmapsto \Ad{s(g)}(h) \,.
\]
Using a different section $s': G\rightarrow H$ to define $\alpha'$ lets us compute
\begin{align*}
\alpha(g,h)\,\alpha'(g,h)^{-1} 
&= s(g)hs(g)^{-1} s'(g) h^{-1} s'(g)^{-1}\\
&= s(g)s(g)^{-1} s'(g) hh^{-1} s'(g)^{-1}
= e
\end{align*}
-- using that $s(g)^{-1}s'(g)$ is in the kernel of $\tau$, thus central. Hence $\alpha$ is independent of the choice of section $s$. Equivariance of $\tau$ is seen using definitions:
\[
\tau \alpha(g,h)  
	= \tau( s(g)hs(g)^{-1} ) 
	= g\, \tau h \, g^{-1}
\]
Noting that elements $h^{-1}s(\tau(h))$ are central in $H$, the Peiffer identity \eqref{peiffer} follows as well,
\[
\alpha(\tau h, h') 
= s(\tau h) h' s(\tau h)^{-1} hh^{-1}
= s(\tau h) s(\tau h)^{-1} hh'h^{-1} 
= hh'h^{-1},
\]
making $(G, H, \tau, \alpha)$ a crossed module as desired. 
The example of the central extension can be seen as the special case of the exact sequence 
\[
1
\longrightarrow
\ker \tau
\longrightarrow
H
\stackrel{\tau}{\longrightarrow}
G
\longrightarrow
\coker \tau
\longrightarrow
1
\]
we get for a crossed module in the case that $\tau$ is an epimorphism. The computations done above highlight that the Peiffer identity forces $\ker \tau$ to be central. 

Taking the other extreme of that sequence with $H$ a normal subgroup and $\tau$ a monomorphism, we get the familiar sequence
\[
1
\longrightarrow
H
\stackrel{\tau}{\longrightarrow}
G
\longrightarrow
G/H
\longrightarrow
1
\]
which gives rise to a crossed module as well with the action $\alpha$ given by conjugation and equivariance as well as the Peiffer identity holding automatically because of that.

%
%

\mysectionstar{References}
\vspace{-1.5em}
\renewcommand{\refname}{}
\printbibliography

\end{document}

%% file: sspbibstyle.tex
\usepackage[style=alphabetic,hyperref=true]{biblatex}
\bibliography{ssp}

\renewcommand\bibfont{\footnotesize\raggedright\frenchspacing}

\DefineBibliographyStrings{english}{
volume = {vol.}
}
\DeclareFieldFormat[article]{title}{\emph{#1}}
\DeclareFieldFormat[incollection]{title}{\emph{#1}}
\DeclareFieldFormat{issn}{}
\DeclareFieldFormat{isbn}{}
\DeclareFieldFormat{url}{}
\DeclareFieldFormat[article]{pages}{\addcolon #1}
\DeclareFieldFormat[book]{pages}{}
\DeclareFieldFormat{booktitle}{#1\isdot}
\DeclareFieldFormat{journaltitle}{#1\isdot}
\DeclareFieldFormat{volumes}{#1~\bibstring{volumes}}
\DeclareFieldFormat{series}{\nopunct#1}

\DeclareFieldFormat{volume}{\bibstring{volume}~#1 of\nopunct}

\renewbibmacro*{issue+date}{%
  \printtext{%
    \iffieldundef{issue}
      {\iffieldundef{month}
         {\printfield{year}}
         {\iffieldundef{day}
            {\printfield{month}%
             \setunit{\addspace}%
             \printfield{year}}
            {\printtext{\bibdate}}}}
      {\printfield{issue}%
       \setunit{\addspace}%
       \printfield{year}}}%
  \newunit}

\newbibmacro*{journal+pages}{%
  \usebibmacro{journal}%
  \setunit*{\addspace}%
  \iffieldundef{series}
    {}
    {\newunit
     \printfield{series}%
     \setunit{\addspace}}%
  \printfield{volume}%
  \printfield[parens]{number}%
  \setunit{\addcolon}%
  \printfield{pages}%
 \setunit{\addcomma\space}%
  \usebibmacro{issue+date}%
  \newunit\newblock
  \usebibmacro{issue}%
  \newunit}
  
\newbibmacro*{publisher+year}{%
  \printlist{publisher}%
  \setunit*{\addcomma\space}%
  \printfield{year}%
  \newunit}

\DeclareBibliographyDriver{article}{%
  \usebibmacro{bibindex}%
  \usebibmacro{author/editor}%
  \setunit{\labelnamepunct}\newblock
  \usebibmacro{title}%
  \newunit
  \usebibmacro{byauthor}%
  \newunit
  \printlist{language}%
  \newunit\newblock
  \usebibmacro{byeditor+others}%
  \newunit\newblock
  \usebibmacro{journal+pages}%
  \newunit\newblock
  \printfield{note}%
  \newunit\newblock
  \printfield{addendum}%
  \newunit\newblock
  \usebibmacro{pageref}%
  \setunit{\addperiod\space}%
  \printfield{eid}%
  \usebibmacro{finentry}%
}

\DeclareBibliographyDriver{online}{%
  \usebibmacro{bibindex}%
  \usebibmacro{author/editor}%
  \setunit{\labelnamepunct}\newblock
  \usebibmacro{title}%
  \setunit{\addcomma\space}%
  \newunit
  \printlist{language}%
  \newunit\newblock
  \usebibmacro{byeditor}%
  \newunit\newblock
  \printfield{version}%
  \newunit
  \printlist{organization}%
 \newunit\newblock  
  \usebibmacro{date}%
  \newunit\newblock
  \usebibmacro{url+urldate}%
 \newunit\newblock
  \printfield{note}%
  \newunit\newblock
  \printfield{addendum}%
  \newunit\newblock
  \usebibmacro{pageref}%
  \setunit{\addperiod\space}%
  \printfield{eid}%
  \usebibmacro{finentry}%
}

\DeclareBibliographyDriver{incollection}{%
  \usebibmacro{bibindex}%
  \usebibmacro{author}%
  \setunit{\labelnamepunct}\newblock
  \usebibmacro{title}%
  \newunit\newblock
  \usebibmacro{in:}%
  \usebibmacro{maintitle+booktitle}%
  \newunit
  \printlist{language}%
  \newunit\newblock
  \usebibmacro{byeditor+others}%
  \newunit\newblock
  \printfield{edition}%
  \newunit
  \iffieldundef{maintitle}
    {\printfield{volume}%
     \printfield{part}}
    {}%
  \newunit
  \printfield{volumes}%
  \newunit\newblock
  \usebibmacro{series+number}%
  \newunit\newblock
  \printfield{note}%
  \newunit\newblock
  \usebibmacro{publisher+year}%
  \addcolon\newblock
  \usebibmacro{chapter+pages}%
  \newunit\newblock
  \printfield{isbn}%
  \newunit\newblock
  \printfield{doi}%
  \newunit\newblock
  \usebibmacro{url+urldate}%
  \newunit\newblock
  \printfield{addendum}%
  \newunit\newblock
  \usebibmacro{pageref}%
  \usebibmacro{finentry}}

\makeatletter
\def\hardspaces@@to#1~#2{%
 #1\ifx#2\@empty\else
   \FROM@SPACE\expandafter\hardspaces@@to
 \fi#2}
\newcommand\spaces@to[2]{%
 \begingroup
   \def\FROM@SPACE{#1}%
   \hardspaces@@to#2~\@empty
 \endgroup}
\renewcommand\mkbibnamefirst[1]{%
 \spaces@to{\,}{#1}}
\makeatother

\DeclareNameFormat{author}{%
  \usebibmacro{name:first-last}{#1}{#4}{#6}{#8}%
  \usebibmacro{name:andothers}}